\documentclass[12pt, reqno, psamsfonts]{amsart}
\pdfoutput=1

\usepackage{amssymb}
\usepackage{amsthm}
\usepackage{amsmath}
\usepackage{latexsym}
\usepackage[T1]{fontenc}
\usepackage[utf8]{inputenc}
\usepackage[russian, french, english]{babel}

\usepackage{graphicx}
\usepackage{wrapfig}
\usepackage[justification=centering, labelfont=bf]{caption}
\usepackage{mathtools}
\usepackage{amsbsy}
\usepackage[inline]{enumitem}
\usepackage{mathrsfs}
\usepackage{array}
\usepackage{multicol}
\usepackage{stmaryrd}
\usepackage{cancel}
\usepackage{lmodern}
\usepackage{mathabx}
\usepackage{upgreek}
\usepackage{titlesec}
\usepackage{titletoc}
\usepackage[spacing=true,kerning=true,babel=true,tracking=true]{microtype}
\usepackage[shortcuts]{extdash}
\usepackage[foot]{amsaddr}
\usepackage[left=1in,right=1in,top=1in,bottom=1in,bindingoffset=0cm]{geometry}
\usepackage{bm}
\usepackage{centernot}
\usepackage{mdframed}
\usepackage[hidelinks]{hyperref}
\usepackage{xspace}
\usepackage[skins]{tcolorbox}
\usepackage{framed}
\usepackage[labelfont=bf,justification=justified]{caption}
\usepackage{tikz}
\usetikzlibrary{shapes,snakes}
\usetikzlibrary{arrows.meta}
\usetikzlibrary{decorations.pathmorphing}
\usetikzlibrary{patterns}
\usepackage{float}
\usepackage{textgreek}

\usepackage[
    sortcites,
    backend=biber, style=alphabetic, sorting=nyt, maxnames=100,backref=true]{biblatex}
\title{\sffamily Sunflowers in Set Systems with Small VC-Dimension}
\date{}

\author{J\'ozsef~Balogh}
\address{\normalfont (JB) Department of Mathematics, University of Illinois at Urbana--Champaign, Urbana, IL, USA. Research supported in part by NSF grants DMS-1764123 and RTG DMS-1937241, FRG DMS-2152488, the Arnold O.~Beckman Research Award (UIUC Campus Research Board RB 24012), and the Simons Fellowship.}
\email{jobal@illinois.edu}

\author{Anton~Bernshteyn}
\address{\normalfont (AB) Department of Mathematics, University of California, Los Angeles, CA, USA. Research supported in part by NSF CAREER grant DMS-2239187 and the Sloan Research Fellowship (2025).}
\email{bernshteyn@math.ucla.edu}

\author{Michelle Delcourt}
\address{\normalfont (MD) Department of Mathematics, Toronto Metropolitan University, Toronto, ON, Canada.  Research supported by NSERC under Discovery Grant No. 2019-042 and the Sloan Research Fellowship (2025).}
\email{mdelcourt@torontomu.ca}

\author{Asaf~Ferber}
\address{\normalfont (AF) Department of Mathematics, University of California, Irvine, CA, USA.  Research is supported in part by NSF grant DMS-1953799, NSF CAREER grant DMS-2146406, the Sloan Research Fellowship (2022), and Air Force grant FA9550-23-1-0298.}
\email{asaff@uci.edu}

\author{Huy~Tuan~Pham}
\address{\normalfont (HP) Department of Mathematics, California Institute of Technology, CA, USA.  Research supported by a Two Sigma Fellowship, a Clay Research Fellowship, and a Stanford Science Fellowship.}
\email{htpham@caltech.edu}

\newtheoremstyle{bfnote}%
{}{}%
{\slshape}{}%
{\bfseries}{\bfseries.}%
{ }%
{\thmname{#1}\thmnumber{ #2}\thmnote{ \ep{\normalfont{}#3}}}

\theoremstyle{bfnote}
\newtheorem{theo}{Theorem}[section]
\newtheorem*{theo*}{Theorem}

\newtheorem{lemma}[theo]{Lemma}
\newtheorem{claim}[theo]{Claim}

\newtheorem{obs}[theo]{Observation}
\newtheorem*{corl*}{Corollary}

\theoremstyle{definition}

\newtheorem*{defn*}{Definition}

\newtheorem*{exmp*}{Example}

\theoremstyle{remark}
\newtheorem*{ques*}{Question}
\newtheorem*{remk*}{Remark}

\newcommand*{\myproofname}{Proof}

\makeatletter
\newcommand{\neutralize}[1]{\expandafter\let\csname c@#1\endcsname\count@}
\makeatother

\newenvironment{theocopy}[1]
{\renewcommand{\thetheo}{\ref{#1}}%
	\neutralize{theo}\phantomsection
	\begin{theo}}
	{\end{theo}}

\newcommand{\0}{\varnothing}
\newcommand{\set}[1]{\{#1\}}
\newcommand{\N}{{\mathbb{N}}}
\newcommand{\Z}{\mathbb{Z}}

\renewcommand{\P}{\mathbb{P}}
\newcommand{\E}{\mathbb{E}}
\renewcommand{\epsilon}{\varepsilon}

\renewcommand{\phi}{\varphi}
\renewcommand{\theta}{\vartheta}
\renewcommand{\leq}{\leqslant}
\renewcommand{\geq}{\geqslant}
\renewcommand{\le}{\leq}

\newcommand{\defeq}{\coloneqq}

\newcommand{\bemph}[1]{{\normalfont#1}} 
\newcommand{\ep}[1]{\bemph{(}#1\bemph{)}} 

\newcommand{\rest}[2]{{{#1}\vert_{#2}}}

\newcommand{\VC}{\mathsf{VC}}
\newcommand{\upset}[1]{\langle #1 \rangle}
\newcommand{\minimal}[1]{{{#1}^\downarrow}}
\newcommand{\core}[2]{\mathsf{core}_{#2}(#1)}

\numberwithin{equation}{section}

\newenvironment{scproof}[1][]{\begin{proof}[\textsc{\upshape{Proof}}#1]}{\end{proof}}

\titleformat{\section}[block]{\large\bfseries\sffamily}{\thesection.}{1ex}{}
\titleformat{\subsection}[block]{\bfseries\sffamily}{\thesubsection.}{1ex}{}
\titleformat{\subsubsection}[runin]{\bfseries}{\bfseries\upshape\sffamily\thesubsubsection.}{1ex}{}

\titlespacing*{\section}{0pt}{*3}{*1}
\titlespacing*{\subsection}{0pt}{*3}{*1}
\titlespacing*{\subsubsection}{0pt}{*2}{*1}

\titlecontents{section}
 [1.5em] %
 {\smallskip}
 {\bfseries\thecontentslabel\hspace{1.02em}}
{\bfseries}
 {\,\,\titlerule*[0.77pc]{}\bfseries\contentspage}
\titlecontents{subsection}
 [4em] %
 {\smallskip}
 {\thecontentslabel\hspace{1.02em}}
{\hspace*{2.32em}}
 {\,\,\titlerule*[0.77pc]{.}\contentspage}

\renewbibmacro{in:}{}

\renewbibmacro*{volume+number+eid}{%
	\printfield{volume}%
	\setunit*{\addnbspace}
	\printfield{number}%
	\setunit{\addcomma\space}%
	\printfield{eid}}

\DeclareFieldFormat[article]{volume}{\textbf{#1}\space}
\DeclareFieldFormat[article]{number}{\mkbibparens{#1}}

\DeclareFieldFormat{journaltitle}{#1,}
\DeclareFieldFormat[thesis]{title}{\mkbibemph{#1}\addperiod}
\DeclareFieldFormat[article, unpublished, thesis]{title}{\mkbibemph{#1},}
\DeclareFieldFormat[book]{title}{\mkbibemph{#1}\addperiod}
\DeclareFieldFormat[unpublished]{howpublished}{#1, }

\DeclareFieldFormat{pages}{#1}

\DeclareFieldFormat[article]{series}{Ser.~#1\addcomma}

\setlist{topsep=3pt,itemsep=3pt}

\begin{document}

    \vspace*{-20pt}

    \maketitle

\begin{abstract}
    A family of $r$ distinct sets $\{A_1,\ldots, A_r\}$ is an {$r$-sunflower} if for all $1 \le i <  j\le r$  and $1 \le i' <  j'\le r$, we have $A_i\cap A_j = A_{i'}\cap A_{j'}$. Erd\H{o}s and Rado conjectured in 1960 that every family $\mathcal{H}$ of $\ell$-element sets of size at least $K(r)^\ell$ contains an $r$-sunflower, where $K(r)$ is some function that depends only on $r$. We prove that if $\mathcal{H}$ is a family of $\ell$-element sets of VC-dimension at most $d$ and $|\mathcal H| > (C r(\log d+\log^\ast\ell))^\ell$ for some absolute constant $C > 0$, then $\mathcal{H}$ contains an $r$-sunflower. This improves a recent result of Fox, Pach, and Suk. When $d=1$, we obtain a sharp bound, namely that $|\mathcal H| > (r-1)^\ell$ is sufficient. Along the way, we establish a strengthening of the Kahn--Kalai conjecture for set families of bounded VC-dimension, which is of independent interest.
\end{abstract}

\section{Introduction}

    All logarithms in this paper are binary and all set systems are finite. A family of $r$ distinct sets $\{A_1,\ldots, A_r\}$ is called an \textbf{$r$-sunflower} (also known as a \textbf{\textDelta-system}) if for all $1 \le i <  j\le r$  and $1 \le i' <  j'\le r$, we have $A_i\cap A_j = A_{i'}\cap A_{j'}$. The intersection $A_i \cap A_j$, common to all pairs $1 \leq i < j \leq r$, is called the \textbf{kernel} of the sunflower. We say that a set family $\mathcal{H}$ is \textbf{$\ell$-bounded} if all sets in $\mathcal{H}$ have size at most $\ell$. Let $f_r(\ell)$ denote the maximum size of an $\ell$-bounded set family $\mathcal{H}$ that does not contain an $r$-sunflower. Erd\H os and Rado~\cite{ER60} proved the so-called \textbf{sunflower lemma}, which states that $f_r(\ell) \leq (r-1)^\ell \, \ell!$, and hence
    \begin{equation}\label{eq:ER_bound}
        f_r(\ell)^{1/\ell} \,=\, O(r \ell).
    \end{equation}
    They also conjectured that $f_r(\ell)^{1/\ell} \leq K(r)$ for some function $K(r)$ that depends only on $r$. 
    In a recent breakthrough,  Alweiss, Lovett, Wu, and Zhang~\cite{ALWZ} improved \eqref{eq:ER_bound} to 
 \[
     f_r(\ell)^{1/\ell} \,=\, O(r^3\log \ell\log\log \ell).
 \]
 The best currently known bound is 
 \begin{equation}\label{eq:BCW_bound}
     f_r(\ell)^{1/\ell} \,=\, O(r \log \ell),
 \end{equation}
 due to Bell, Chueluecha, and Warnke \cite{BCW21} (see also \cite{Rao20} for an intermediate result by Rao). Perhaps even more importantly than making progress on the Erd\H{o}s--Rado sunflower conjecture, Alweiss, Lovett, Wu, and Zhang \cite{ALWZ} invented the so-called \textbf{spread method}, which was further developed by Frankston, Kahn, Narayanan, and Park \cite{FKNP21} and inspired the proof of the Kahn--Kalai conjecture \cite{KK07} and Talagrand's selector process conjecture \cite{Tal1, Tal2, Tal3} by Park and Pham \cite{PP24, PP24a}. 
 These developments have led to a number of other breakthroughs in probabilistic and extremal combinatorics \cite{pham2025sharp,
spread1,spread2,spread3,spread4,spread6,spread7,spread8,spread9,spread10,spread11,spread12,spread13,spread14,kelly2024optimal,nenadov2024spread,allen2025robustness,han2025perturbation}. 
We direct the reader to the surveys \cite{Jsurvey} by Park and \cite{Wsurvey} by Perkins and to Pham's PhD thesis \cite{PhamThesis} for an overview of the area. 

    The \textbf{Vapnik--Chervonenkis dimension}, or \textbf{VC-dimension} for short, of a set family $\mathcal{H}$ on a ground set $X$ is the maximum $d$ for which there exists a subset $T \subseteq X$ of size $d$ such that for all $A \subseteq T$, there is some $S \in \mathcal{H}$ with $T \cap S = A$. 
    We denote the VC-dimension of a set family $\mathcal{H}$ by $\VC(\mathcal{H})$. This notion was developed by Vapnik and Chervonenkis \cite{VC71} and plays an important role in statistical learning theory \cite{learning1,learning2}, discrete and computational geometry \cite{geom1,geom2,geom3}, and model theory \cite{NIP}. In recent years, it has also found applications in extremal combinatorics, wherein a bound on VC-dimension is used to obtain improved extremal results; see \cite{extremal1,extremal2,extremal3,extremal4,extremal5,extremal6,extremal7,extremal8} for a selection of papers belonging to this research stream. One such application is an improved bound on the sunflower problem for set families of bounded VC-dimension, recently established by Fox, Pach, and Suk~\cite{FPS23}: 
    



    \begin{theo}[{Fox--Pach--Suk \cite[Theorem 1.3]{FPS23}}]\label{FPSthm}
        The following holds with $C = 2^{10}$. Let $\mathcal{H}$ be an $\ell$-bounded set system with $\VC(\mathcal{H}) \leq d$, where $d \geq 2$. 
        If
        \[
            |\mathcal{H}|^{1/\ell} \,\geq\,C^{(dr)^{2 \log^\ast \ell}},
        \]
        then $\mathcal{H}$ contains an $r$-sunflower.
    \end{theo}

    Here $\log^\ast \ell$ is the \textbf{iterated logarithm} of $\ell$, i.e., the number of times the binary logarithm function needs to be applied to $\ell$ before the result becomes at most $1$. For constant $d$ and $r$, the bound on $|\mathcal{H}|^{1/\ell}$ given by Theorem~\ref{FPSthm} is doubly-exponential in $\log^\ast \ell$. In our main result we replace this doubly-exponential dependence 
    by a linear one: 

    \begin{theo}[Main result: a sunflower lemma for set families of bounded VC-dimension]\label{theo:sunflower}
        There exists a universal constant $C > 0$ such that the following holds. Let $\mathcal{H}$ be an $\ell$-bounded set system with $\VC(\mathcal{H}) \leq d$, where $d \geq 2$. If
        \[
            |\mathcal{H}|^{1/\ell} \,\geq\, C r (\log d + \log^\ast \ell),
        \]
        then $\mathcal{H}$ contains an $r$-sunflower.
    \end{theo}

    Notice that our result also improves Theorem~\ref{FPSthm} in its dependence on $r$ and $d$. In particular, since every $\ell$-bounded set family $\mathcal{H}$ satisfies $\VC(\mathcal{H}) \leq \ell$, the bound given by our theorem is always at least as good as the general result \eqref{eq:BCW_bound}, up to a constant factor. Furthermore, for $d = \ell^{o(1)}$, we get a strict asymptotic improvement on \eqref{eq:BCW_bound}. By contrast, the bound in Theorem~\ref{FPSthm} can only be better than \eqref{eq:BCW_bound} when $d \ll \log \log \ell$.
A preliminary version of Theorem~\ref{theo:sunflower} with a looser (namely, exponential) dependence on $d$ appeared previously in the fifth author's PhD thesis \cite[\S4]{PhamThesis}.



    In the special case $\VC(\mathcal{H}) = 1$ we can be even more precise. 
In~\cite{FPS23} the following construction was described to show that there exists an $\ell$-bounded set family of size $(r-1)^{\ell}$ with VC-dimension $1$ and no $r$-sunflower. Fix a rooted complete $(r-1)$-ary tree $T$ with the root on level $0$ and with $(r-1)^{\ell}$
leaves on level $\ell$. Let $\mathcal{H}$ be the family of the edge sets of the leaf-to-root paths in $T$. It is straightforward to check that $\mathcal{H}$ has the required properties. Fox, Pach, and Suk complemented this construction by the following bound:

    \begin{theo}[{Fox--Pach--Suk \cite[Theorem 1.2]{FPS23}}]\label{FPSthm1}
        Let $\mathcal{H}$ be an $\ell$-bounded set system with $\VC(\mathcal{H}) \leq 1$. If $|\mathcal{H}| > r^{10\ell}$, then $\mathcal{H}$ contains an $r$-sunflower.
    \end{theo}

%


 In our second result, we improve this by showing that the tree construction is optimal:

\begin{theo}[Sharp bound for 1-dimensional families]\label{pd1r3}
    Let $\mathcal{H}$ be an $\ell$-bounded set system with $\VC(\mathcal{H}) \leq 1$. If $|\mathcal{H}| > (r-1)^{\ell}$, then $\mathcal{H}$ contains an $r$-sunflower.
\end{theo}

    
    Let us now say a few words about the machinery employed in the proof of Theorem~\ref{theo:sunflower}. To begin with, we need some notation. For a set system $\mathcal{H}$ on a ground set $X$, define
    \[
        \upset{\mathcal{H}} \,\defeq\, \set{A \subseteq X \,:\,  \exists S \in \mathcal{H} \, \text{ such that } A \supseteq S}.
    \]
    For $p \in [0,1]$, let $X_p$ be the probability distribution on $\mathcal{P}(X)$, the power set of $X$, where a random subset of $X$ is generated by including each element independently with probability $p$. For convenience, if $p > 1$, we let $X_p \defeq X_1$ (i.e., if $W \sim X_p$, then $\P[W = X] = 1$). The following result (without explicit constants) was conjectured by Kahn and Kalai \cite{KK07}. Their conjecture was recently verified by Park and Pham \cite{PP24} (a weaker, fractional variant was established earlier by Frankston, Kahn,  Narayanan, and Park~\cite{FKNP21}). The version we state below is due to Bell \cite{B23} and includes a refinement of the dependence on $\epsilon$.

    \begin{theo}[{Kahn--Kalai conjecture~\cite{KK07}; Park--Pham \cite{PP24}, Bell \cite{B23}}]\label{theo:KKBell}
        Let $\mathcal{H}$ be an $\ell$-bounded set system on a ground set $X$ and let $q \in [0,1]$. Then at least one of the following statements holds.
        \begin{enumerate}[label=\ep{$\text{\normalfont{\small{KK}}}_{\arabic*}$}]
            \item\label{item:KK1} There is a set system $\mathcal{F}$ with $\mathcal{H} \subseteq \upset{\mathcal{F}}$ and $\displaystyle \sum_{F \in \mathcal{F}} q^{|F|} \leq \frac{1}{2}$.

            \item\label{item:KK2} For all $\epsilon \in (0, \frac{1}{2}]$, if $p = 48 q \log(\ell/\epsilon)$ and $W \sim X_p$, then $\P[W \in \upset{\mathcal{H}}] > 1 - \epsilon$.
        \end{enumerate}
    \end{theo}

   An important general problem is to determine when the logarithmic factor in the definition of $p$ in \ref{item:KK2} can be removed; see, e.g., \cite{spread3,spread4,spread14}. We show that for families $\mathcal{H}$ of bounded VC-dimension, the logarithmic dependence on $\ell$ can be reduced to just $\log^* \ell$: 

    \begin{theo}[Kahn--Kalai for set families of bounded VC-dimension]\label{theo:KK}
        There is a universal constant $A > 0$ with the following property. Let $\mathcal{H}$ be an $\ell$-bounded set system on a ground set $X$ with $\VC(\mathcal{H}) \leq d$ and let $q \in [0,1]$. Then at least one of the following statements holds.
        \begin{enumerate}[label=\ep{$\text{\normalfont{\small{VC}}}_{\arabic*}$}]
            \item\label{item:VC1} There is a set system $\mathcal{F}$ with $\mathcal{H} \subseteq \upset{\mathcal{F}}$ and $\displaystyle \sum_{F \in \mathcal{F}} q^{|F|} \leq \frac{2}{3}$.

            \item\label{item:VC2} For all $\epsilon \in (0, \frac{1}{2}]$, if $p = A q (\log(d/\epsilon) + \log^\ast \ell)$ and $W \sim X_p$, then $\P[W \in \upset{\mathcal{H}}] > 1 - \epsilon$.
        \end{enumerate}
    \end{theo}
    
    For families of VC-dimension $1$, Theorem~\ref{theo:KK} can be further improved by completely eliminating the dependence on $\ell$:

    \begin{theo}[Kahn--Kalai for set families of VC-dimension 1]\label{theo:KKdim1}
        There is a universal constant $A > 0$ with the following property. Let $\mathcal{H}$ be a set system on a ground set $X$ with $\VC(\mathcal{H}) \leq 1$ and let $q \in [0,1]$. Then at least one of the following statements holds.
        \begin{enumerate}[label=\ep{$\text{\normalfont{\small{1d}}}_{\arabic*}$}]
            \item\label{item:VC1dim1} There is a set system $\mathcal{F}$ with $\mathcal{H} \subseteq \upset{\mathcal{F}}$ and $\displaystyle \sum_{F \in \mathcal{F}} q^{|F|} \leq \frac{2}{3}$.

            \item\label{item:VC2dim1} For all $\epsilon \in (0, \frac{1}{2}]$, if $p = A q \log(1/\epsilon)$ and $W \sim X_p$, then $\P[W \in \upset{\mathcal{H}}] > 1 - \epsilon$.
        \end{enumerate}
    \end{theo}

    The rest of the paper is organized as follows. Theorem~\ref{pd1r3} is proved in \S\ref{sd1r3}. Then, in \S\ref{sp13}, we prove Theorem~\ref{theo:KK}. Theorem~\ref{theo:sunflower} follows from Theorem~\ref{theo:KK} in a fairly straightforward manner, inspired by the proof of \eqref{eq:BCW_bound} due to Bell, Chueluecha, and Warnke \cite{BCW21}. We present this argument in \S\ref{stheo}. Finally, in \S\ref{sd1genr} we show how to modify the proof of Theorem~\ref{theo:KK} in the case $d = 1$ to establish Theorem~\ref{theo:KKdim1}. 

    \medskip

    \noindent \textbf{Acknowledgments.} The first four authors sincerely thank the American Institute of Mathematics for hosting the SQuaREs program ``Containers from Different Angles,'' where most of their work on this project was done. They also thank Caroline Terry and  Anush~Tserunyan for useful discussions on the project. The first author is grateful to Robert Krueger for fruitful discussions on the Kahn--Kalai conjecture/Park--Pham theorem.  The third author would like to thank Luke Postle and Tom Kelly for a number of enlightening conversations on related topics. The fourth author is grateful to Marcelo Sales for fruitful discussions. Finally, we are thankful to the anonymous referees for carefully reading the manuscript and providing helpful comments.

    \section{Sunflowers in 1-dimensional families: Proof of Theorem~\ref{pd1r3}}\label{sd1r3}

    The main part of the proof of Theorem~\ref{pd1r3} is the following lemma:
    
    \begin{lemma}\label{struct}
        Let $\mathcal{H}$ be an $\ell$-bounded family of sets on a ground set $X$ with $|\mathcal{H}| > (r-1)^\ell$. If $\mathcal{H}$ contains no $r$-sunflower, then there exist distinct elements $x$, $y \in X$ and sets $S_x$, $S_y$, $S_{xy} \in \mathcal{H}$ such that $\set{x,y} \cap S_x = \set{x}$, $\set{x,y} \cap S_y = \set{y}$, and $\set{x,y} \cap S_{xy} = \set{x,y}$.
    \end{lemma}
    \begin{scproof} We proceed by induction on $\ell$. The base case $\ell=1$ holds vacuously as every $1$-bounded set family with at least $r$ members contains an $r$-sunflower with an empty kernel. Now suppose $\ell \geq 2$ and the statement holds when $\ell$ is replaced by $\ell - 1$. Without loss of generality, we assume that every set in $\mathcal{H}$ has exactly $\ell$ elements (otherwise we may add $|S|-\ell$ unique new elements to each set $S \in \mathcal{H}$). In particular, the sets in $\mathcal{H}$ are nonempty.

    Let $F_1$, \ldots, $F_t \in \mathcal{H}$ be a maximal sequence of pairwise disjoint sets in $\mathcal{H}$ (it is possible that $t = 1$). Note that $t \leq r-1$, since otherwise $\set{F_1, \ldots, F_t}$ would be a sunflower (with an empty kernel) of size at least $r$. By the choice of $F_1$, \ldots, $F_t$, every set $F \in \mathcal{H}$ intersects at least one of these sets. Therefore, there exists an index $i^* \in [t]$ such that the family $\mathcal{H}_{i^*} \defeq \{F\in \mathcal{H} \,:\, F\cap F_{i^*} \neq \0\}$ satisfies $|\mathcal{H}_{i^*}| \geq |\mathcal{H}|/t > (r-1)^{\ell - 1}$.

    If there exist $F$, $F'\in \mathcal{H}_{i^*}$ such that the sets $E \defeq F\cap F_{i^*}$ and $E'\defeq F'\cap F_{i^*}$ are not comparable with respect to inclusion, then choosing any $x\in E\setminus E'$ and $y\in E'\setminus E$ and letting $S_x \defeq F$, $S_y \defeq F'$, and $S_{xy} \defeq F_{i^*}$ completes the proof. Therefore, we may assume that the family of intersections $\mathcal{E}\defeq\{F\cap F_{i^*}\,:\, F \in \mathcal{H}_{i^*}\}$ forms a chain under inclusion. 

    Let $E\in \mathcal{E}$ be the minimal element in the chain $\mathcal{E}$. Then $\0 \neq E \subseteq F$ for all $F \in \mathcal{H}_{i^*}$. Consider the set system $\mathcal{H}' \defeq \set{F \setminus E \,:\, F \in \mathcal{H}_{i^*}}$. Since $E \neq \0$, $\mathcal{H}'$ is $(\ell - 1)$-bounded and satisfies $|\mathcal{H}'| = |\mathcal{H}_{i^*}| > (r-1)^{\ell - 1}$. Moreover, it does not contain an $r$-sunflower, since an $r$-sunflower in $\mathcal{H}'$ corresponds to an $r$-sunflower in $\mathcal{H}$ (with a kernel including $E$). Therefore, by the inductive hypothesis, the conclusion of the lemma holds for $\mathcal{H}'$, which implies that it also holds for $\mathcal{H}$. 
\end{scproof}

    We are now ready to prove Theorem~\ref{pd1r3}, which we restate here for convenience.

    \begin{theocopy}{pd1r3}
        Let $\mathcal{H}$ be an $\ell$-bounded set system on a ground set $X$ with $\VC(\mathcal{H}) \leq 1$. If $|\mathcal{H}| > (r-1)^{\ell}$, then $\mathcal{H}$ contains an $r$-sunflower.
    \end{theocopy}
    \begin{scproof}
        We may assume that $r \geq 3$, as for $r \in \set{1,2}$ the theorem holds trivially. The proof is by induction on $\ell$. If $\ell = 1$, then $\mathcal{H}$ contains an $r$-sunflower with an empty kernel. Now suppose that $\ell \geq 2$ and the theorem holds when $\ell$ is replaced by $\ell - 1$.

        Suppose 
        that $\mathcal{H}$ contains no $r$-sunflower. By Lemma~\ref{struct}, there exist elements $x$, $y \in X$ and sets $S_x$, $S_y$, $S_{xy} \in \mathcal{H}$ such that $\set{x,y} \cap S_x = \set{x}$, $\set{x,y} \cap S_y = \set{y}$, and $\set{x,y} \cap S_{xy} = \set{x,y}$. Since $\VC(\mathcal{H}) \leq 1$, there must be some $A \subseteq \set{x,y}$ such that $A \neq \set{x,y} \cap S$ for all $S \in \mathcal{H}$. Due to the existence of the sets $S_x$, $S_y$, $S_{xy}$, the only option for $A$ is $A = \0$. In other words, every set $S \in \mathcal{H}$ has a nonempty intersection with $\set{x,y}$. It follows that there exists $z \in \{x,y\}$ such that the subfamily $\mathcal{H}_z\defeq\left\{S\in \mathcal{H} \,:\, z\in S\right\}$
satisfies  
    $
        |\mathcal{H}_z| \geq |\mathcal{H}|/2 > (r-1)^{\ell-1}
    $ (here we use our assumption that $r \geq 3$).
    
    Consider the $(\ell-1)$-bounded family $\mathcal{H}' \defeq \set{S \setminus \set{z} \,:\, S \in \mathcal{H}_{z}}$. It is clear that $\VC(\mathcal{H}') \leq \VC(\mathcal{H}) \leq 1$. As $|\mathcal{H}'| = |\mathcal{H}_{z}| > (r-1)^{\ell - 1}$, we conclude that $\mathcal{H}'$ contains an $r$-sunflower by the inductive hypothesis. But an $r$-sunflower in $\mathcal{H}'$ corresponds to an $r$-sunflower in $\mathcal{H}$ (with a kernel containing $z$), so $\mathcal{H}$ has an $r$-sunflower as well, and we are done.
    \end{scproof}

\section{Kahn--Kalai for low-dimensional families: Proof of Theorem~\ref{theo:KK}}\label{sp13}

    \subsection{Preliminaries}\label{subsec:prelim}

    Throughout the proof of Theorem~\ref{theo:KK}, we will frequently invoke the following simple observation, which follows immediately from the definition of VC-dimension:

\begin{obs}\label{obs:trace}
 Let $\mathcal{H}$ be a set system on a ground set $X$ and let $U \subseteq X$. Define the {\upshape\textbf{trace}} of $\mathcal{H}$ on $U$ by $\rest{\mathcal{H}}{U} \defeq \{S \cap U:\ S\in\mathcal{H}\}$. Then $\VC(\rest{\mathcal{H}}{U}) \leq \VC(\mathcal{H})$.
\end{obs}

    The following famous result of Sauer~\cite{Sa72}, Perles and Shelah~\cite{Sh72}, and  Vapnik and Chervonenkis~\cite{VC71} bounds the cardinality of a set family in terms of its VC-dimension:

 \begin{lemma}[{Sauer--Shelah lemma \cite{Sa72,Sh72,VC71}; see \cite[Lemma~6.10]{learning2}}]\label{lemma:SS}
   Let $\mathcal{H}$ be a set system on a ground set $X$ with $\VC(\mathcal{H}) \leq d \leq |X|$. Then
   \[
    |\mathcal{H}| \,\le\, \sum_{i=0}^d\binom{|X|}{i} \,\le\, \left(\frac{e|X|}{d}\right)^d.
    \]
 \end{lemma}

    Importantly, due to Observation~\ref{obs:trace}, the bound given by Lemma~\ref{lemma:SS} also applies to the trace of $\mathcal{H}$ on any subset of $X$. This fact is the only property of VC-dimension we shall rely on in our proof of Theorem~\ref{theo:KK}.
   
    Before moving on, we make a technical remark regarding the $\log^\ast$ function. 
    Our proof of Theorem~\ref{theo:KK} will proceed by induction on $\ell$. As a result, it is somewhat inconvenient for our purposes that $\log^\ast$ is piecewise constant, as we would like to use the inequality $\log^\ast \ell' < \log^\ast \ell$ for the value $\ell' < \ell$ to which the inductive hypothesis is applied. To solve this issue, we slightly alter the definition of $\log^\ast$ to  ``smooth it out,'' as follows. The usual definition is that $\log^\ast\ell = t$ if $\ell$ is between the towers of $2$s of heights $t-1$ and $t$. For example, $\log^\ast\ell = 4$ for 
    \[
        16 \,=\, 2^{2^{2}} \,<\, \ell \,\leq\, 2^{2^{2^2}} \,=\,65536.
    \]
    We modify this by splitting the interval between the two towers of $2$s at the point where the top $2$ in the first tower is replaced by a $3$ and letting $\log^\ast \ell \defeq t$ in the first half of the interval and $\log^\ast \ell \defeq t + 1/2$ in the second half. For example, we have 
    \[
        \log^\ast\ell = 4 \ \Longleftrightarrow \ 2^{2^{2}} \,<\, \ell \,\leq\, 2^{2^{3}} \quad \text{and} \quad \log^\ast \ell = 4.5 \ \Longleftrightarrow \ 2^{2^{3}} \,<\, \ell \,\leq\, 2^{2^{2^2}}.
    \]
With this convention, the following holds (we omit the routine proof):
\begin{claim}\label{logrec}
For $x > 8$, 
the following relation holds:
  $\displaystyle \log^\ast(x^2)  + \frac{1}{2}  \le   \log^\ast(2^x)$. 
 \end{claim}

    Clearly, since we may adjust the values of the constant factors in our results, it does not matter whether we use our definition of $\log^\ast$ or the standard one.

    \subsection{A double counting argument}\label{subsec:double}

    The heart of the proof of Theorem~\ref{theo:KK} is in a certain double counting argument, analogous to the one employed by Park and Pham in their proof of the Kahn--Kalai conjecture \cite[Lemma 2.1]{PP24}, but using the VC-dimension of $\mathcal{H}$ to strengthen the bound. 
    
    Let $\mathcal{H}$ be a set system on a ground set $X$. For a set system $\mathcal{S}$, we define $\minimal{\mathcal{S}}$ to be  the set of all inclusion-minimal members of $\mathcal{S}$. Given a subset $W \subseteq X$, we let
    \[
        \mathcal{H}_W \,\defeq\, \minimal{\set{S \setminus W \,:\, S \in \mathcal{H}}} \,=\, \minimal{(\rest{\mathcal{H}}{X \setminus W})},
    \]
    where $\rest{\mathcal{H}}{X \setminus W}$ is the trace of $\mathcal{H}$ on $X \setminus W$ (as defined in Observation~\ref{obs:trace}).
    For a set $A \in \upset{\mathcal{H}}$, the \textbf{$\mathcal{H}$-core} of $A$, denoted by $\core{A}{\mathcal{H}}$, is the intersection of all sets $S \in \mathcal{H}$ such that $S \subseteq A$ (note that there exists at least one such set $S$ because $A \in \upset{\mathcal{H}}$). 

    \begin{figure}[t]
			\centering
			\begin{tikzpicture}

            \begin{scope}
                    \clip (0,1) rectangle (-3,-1);
                    \filldraw[thick,fill=gray!70] (0,0) ellipse (1cm and 0.5cm);
                \end{scope}

                \begin{scope}
                \clip (0,-1) rectangle (3,1);

                \begin{scope}
                    \clip[yslant=0.4] (0,0) ellipse (2cm and 0.5cm);
                    \fill[gray!70,yslant=-0.4] (0,0) ellipse (2cm and 0.5cm);
                \end{scope}

                \end{scope}

            \draw[pattern={north east lines},pattern color=gray!40] (0,-2) rectangle (3,2);
            \draw[thick] (0,0.5) -- (0,-0.5);

                \begin{scope}
                \clip (0,-1) rectangle (3,1);

                \draw (0,0) ellipse (2cm and 0.5cm);
                \draw[yslant=0.4] (0,0) ellipse (2cm and 0.5cm);
                \draw[yslant=-0.4] (0,0) ellipse (2cm and 0.5cm);

                \end{scope}

                \draw (0,-2) -- (-2,-2) -- (-2,2) -- (0,2);
                \node at (-0.5,0) {\small$F$};
                \node at (2.5,-1.6) {\small$W$};
                \node at (2.3,0.8) {\small$S_1$};
                \node at (2.3,0) {\small$S_2$};
                \node at (2.3,-0.8) {\small$S_3$};
                \node at (-1.3,-1.6) {\small$X\setminus W$};

			\end{tikzpicture}
        \caption{An illustration for Lemma~\ref{lemma:core}. Here $F \in \mathcal{H}_W$, the sets $S_1$, $S_2$, $S_3$ are the members of $\mathcal{H}$ contained in $W \cup F$ (note that each of them includes $F$), and the gray region is $\core{W \cup F}{\mathcal{H}}$.
        }\label{fig:core}
	\end{figure}

    \begin{lemma}\label{lemma:core}
        For any $W \subseteq X$ and $F \in \mathcal{H}_W$, we have  $W \cup F \in \upset{\mathcal{H}}$ and $F \subseteq  \core{W \cup F}{\mathcal{H}}$.
    \end{lemma}
    \begin{scproof}
        Since $F \in \mathcal{H}_W$, there is some $S \in \mathcal{H}$ such that $F = S \setminus W$, from which it follows that $S \subseteq W \cup F$ and hence $W \cup F \in \upset{\mathcal{H}}$. The minimality of $F$ implies that $F$ is contained in every set $S \in \mathcal{H}$ with $S \subseteq W \cup F$, therefore $F \subseteq \core{W \cup F}{\mathcal{H}}$. See Fig.~\ref{fig:core} for an illustration.
    \end{scproof}

    Given a set $W \subseteq X$, for each $S \in \mathcal{H}$, we let $F_{\mathcal{H},W}(S)$ be an arbitrary fixed set $F \in \mathcal{H}_W$ such that $F \subseteq S \setminus W$, which exists by the definition of $\mathcal{H}_W$. The sets $F_{\mathcal{H},W}(S)$ are called {minimum fragments} in \cite{PP24} and are central to the proof of the Kahn--Kalai conjecture. To improve the bound using the VC-dimension of $\mathcal{H}$, we shall also utilize the following larger sets:
    \[
        F^\ast_{\mathcal{H},W}(S) \,\defeq\, S \cap \core{W \cup F_{\mathcal{H},W}(S)}{\mathcal{H}}.
    \]
    It follows from the definitions and Lemma~\ref{lemma:core} that
    \begin{equation}\label{eq:inclusions}
        F_{\mathcal{H},W}(S) \,\subseteq\, F^\ast_{\mathcal{H},W}(S) \,\subseteq\, S \qquad \text{and} \qquad F^\ast_{\mathcal{H},W}(S) \setminus W = F_{\mathcal{H},W}(S).
    \end{equation}
    For  $t \geq 0$ we define
 \begin{align*}
        \mathcal{H}_{W, t}^\mathsf{small} \,&\defeq\, \set{F_{\mathcal{H},W}(S) \,:\, S \in \mathcal{H},\,|F_{\mathcal{H},W}(S)| < t} \\
        \mathcal{H}_{W, t}^\mathsf{large} \,&\defeq\, \set{F^\ast_{\mathcal{H},W}(S) \,:\, S \in \mathcal{H}, \, |F_{\mathcal{H},W}(S)| \geq t}.
    \end{align*}
    The following lemma records some basic properties of $\mathcal{H}_{W, t}^\mathsf{small}$ and $\mathcal{H}_{W, t}^\mathsf{large}$:

    \begin{lemma}\label{property}
        For all $W\subseteq X$ and $t \geq 0$, the following holds.


        \begin{enumerate}[label=\ep{\normalfont\arabic*}]
            \item\label{item:undercover} $\mathcal{H} \subseteq \upset{\mathcal{H}_{W,t}^\mathsf{small} \cup \mathcal{H}_{W,t}^\mathsf{large}}$.


            \item\label{item:small} The set system $\mathcal{H}_{W,t}^\mathsf{small}$ is $t$-bounded.


            \item\label{item:induction} If for some set $W' \subseteq X$, we have $W' \in \upset{\mathcal{H}_{W,t}^\mathsf{small}}$, then $W \cup W' \in \upset{\mathcal{H}}$.


            \item\label{item:VC} $\VC(\mathcal{H}_{W,t}^\mathsf{small}) \leq \VC(\mathcal{H})$. 

            \item\label{item:large_is_large} For all $F \in \mathcal{H}_{W,t}^\mathsf{large}$, we have $|F \setminus W| \geq t$.
        \end{enumerate}
    \end{lemma}
    \begin{scproof}
        Follows immediately from the definitions, relations \eqref{eq:inclusions}, and Observation~\ref{obs:trace}.
    \end{scproof}

    We are now ready to state our main double counting lemma:
    
    \begin{lemma}\label{lemma:count}
        Suppose the VC-dimension of $\mathcal{H}$ is at most $d \leq \ell$. Let $p$, $q \in [0,1]$ be such that $p \geq 2q$ and let $W \sim X_p$ be a random subset of $X$. Then, for all $t \geq 0$,
        \[
            \E\left[\sum_{F \,\in\, \mathcal{H}_{W,t}^{\mathsf{large}}} q^{|F|}\right] \,\leq\, 2 \, \left(\frac{e\ell}{d}\right)^d \, \left(\frac{q}{p}\right)^{t}.
        \]
    \end{lemma}
    \begin{scproof}
        Let $n \defeq |X|$ 
        and write
        \begin{align*}
            \E\left[\sum_{F \,\in\, \mathcal{H}_{W,t}^{\mathsf{large}}} q^{|F|}\right] \,&=\, \sum_{W \subseteq X} \sum_{F \,\in\, \mathcal{H}_{W,t}^{\mathsf{large}}} q^{|F|} \, p^{|W|} \, (1 - p)^{n - |W|} \\
            &=\, \sum_{w = 0}^{n} \sum_{k = t}^\ell q^k p^w (1-p)^{n - w} \cdot \left|\left\{(W, F) \,:\, |W| = w, \, F \in \mathcal{H}_{W,t}^{\mathsf{large}}, \, |F| = k\right\}\right|.
        \end{align*}
        Fix $k \in [t, \ell]$ and let us count the pairs $(W, F)$ as above (this step is the crux of the proof). Say $|F \setminus W| = j$, so $t \leq j \leq k$ by Lemma~\ref{property}\ref{item:large_is_large}. Following Park and Pham's proof of \cite[Lemma~2.1]{PP24}, we first fix the union $Z \defeq W \cup F$, for which there are at most ${n \choose w + j}$ choices. Importantly, the $\mathcal{H}$-core $\core{Z}{\mathcal{H}}$ is uniquely determined by $Z$ and has size at most $\ell$. Since $F \in \mathcal{H}_{W,t}^{\mathsf{large}}$, we have $F = S \cap \core{Z}{\mathcal{H}}$ for some $S \in \mathcal{H}$ by definition. Thus, as $\VC(\mathcal{H}) \leq d \leq \ell$, the Sauer--Shelah lemma (Lemma~\ref{lemma:SS}) applied to $\rest{\mathcal{H}}{\core{Z}{\mathcal{H}}}$ in place of $\mathcal{H}$ implies that there are at most $(e\ell/d)^d$ choices for $F$ given $Z$. Finally, given $F$ with $|F| = k$, there are at most ${k \choose j}$ choices for the subset $F \setminus W \subseteq F$ of size $j$, which, together with $Z$, determines $W$. Putting everything together, the quantity we want to bound is at most
        \begin{align}
            &\sum_{w = 0}^{n} \sum_{k = t}^\ell q^k p^w (1-p)^{n - w} \sum_{j = t}^k {n \choose w + j} \,\left(\frac{e\ell}{d}\right)^d\, {k \choose j} \nonumber\\
            =\, &\left(\frac{e\ell}{d}\right)^d \sum_{w =0}^{n} p^w (1-p)^{n-w} \sum_{j = t}^\ell {n \choose w + j} \sum_{k = j}^\ell q^k \,{k \choose j}. \label{eq:1}
        \end{align}
        Next we notice that
        \[
            \sum_{k = j}^\infty q^k \,{k \choose j} \,=\, \frac{q^j}{(1-q)^{j+1}},
        \]
        so the last expression in \eqref{eq:1} is bounded above by
        \begin{align}
            &\left(\frac{e\ell}{d}\right)^d \sum_{w =0}^{n} p^w (1-p)^{n-w} \sum_{j = t}^\ell {n \choose w + j} \frac{q^j}{(1-q)^{j+1}} \nonumber\\
            =\, &\left(\frac{e\ell}{d}\right)^d \sum_{j = t}^\ell \frac{q^j}{(1-q)^{j+1}} \sum_{w = 0}^{n} {n \choose w + j} p^w (1-p)^{n - w} \nonumber\\
            =\, &\left(\frac{e\ell}{d}\right)^d \sum_{j = t}^\ell \frac{q^j(1-p)^j}{(1-q)^{j+1}p^j} \sum_{w = 0}^{n} {n \choose w + j} p^{w+j} (1-p)^{n - w - j}. \label{eq:2}
        \end{align}
        Observe that
        \[
            \sum_{w = 0}^{n} {n \choose w + j} p^{w+j} (1-p)^{n - w - j} \,\leq\, \sum_{s \in \Z} {n \choose s} p^s (1-p)^{n-s} \,=\, 1,
        \]
        so the last expression in \eqref{eq:2} is at most
        \begin{align*}
            \left(\frac{e\ell}{d}\right)^d \sum_{j = t}^\ell \frac{q^j(1-p)^j}{(1-q)^{j+1}p^j} \,\leq\, \left(\frac{e\ell}{d}\right)^d \frac{1}{1-q}\sum_{j = t}^\infty \left(\frac{q(1-p)}{(1-q)p}\right)^j \,=\, \left(\frac{e\ell}{d}\right)^d \frac{p}{p-q} \left(\frac{q(1-p)}{(1-q)p}\right)^t.
        \end{align*}
        Finally, note that since $p \geq 2q$ by assumption, we have $p/(p-q) \leq 2$ and $(1-p)/(1-q) \leq 1$, so the last expression is bounded above by the desired quantity
        \[
            2 \,\left(\frac{e\ell}{d}\right)^d\, \left(\frac{q}{p}\right)^t. \qedhere
        \]
    \end{scproof}

    \subsection{Induction on $\ell$}\label{subsec:proof_of_KK}

    We now have all the ingredients necessary to prove Theorem~\ref{theo:KK}. We first establish a slightly more technical statement designed to facilitate the inductive argument. For $d$, $\ell \in \N$, let
    \[
        \lambda_d(\ell) \,\defeq\, \begin{cases}
            \log^\ast \ell + 2 &\text{if } \ell > 2^{3d},\\
            \log^\ast \ell + 1 &\text{if } 9d^2 < \ell \leq 2^{3d}, \\
            \log^\ast \ell &\text{if }  \ell \leq 9d^2.
        \end{cases}
    \]
    Recall that here we use the variant of the $\log^\ast$ function described in \S\ref{subsec:prelim}.

    \begin{theo}\label{theo:KKvariant}
        There is a universal constant $C > 0$ with the following property. Let $\mathcal{H}$ be an $\ell$-bounded set system on a ground set $X$ with $\VC(\mathcal{H}) \leq d$ 
        and let $q \in [0,1]$. Choose $\epsilon \in (0, \frac{1}{2}]$ and let $\ell_0 \defeq 300\, (d/\epsilon)^3$. Then at least one of the following statements holds.
        \begin{enumerate}[label=\ep{$\text{\normalfont{\small{Ind}}}_{\arabic*}$}]
            \item\label{item:VC1ind} There is a set system $\mathcal{F}$ with
            \[
                \mathcal{H} \,\subseteq\, \upset{\mathcal{F}} \quad \text{and} \quad \sum_{F \in \mathcal{F}} q^{|F|} \,\leq\, \frac{1}{2} + 2 \sum_{i = \ell_0}^\ell \left(\frac{e}{di}\right)^d.
            \]

            \item\label{item:VC2ind} For $p = C q (\log (d/\epsilon) + \lambda_d(\ell))$ and $W \sim X_p$, we have 
            \[
                \P[W \in \upset{\mathcal{H}}] \,>\, 1 - \epsilon - \sum_{i = \ell_0}^\ell i^{-d}.
            \]
        \end{enumerate}
    \end{theo}
    \begin{scproof}
        Let $C > 0$ be a constant that we will assume to be sufficiently large for all of the inequalities that appear in the sequel to hold. We fix $d$ and prove the result by induction on $\ell$. Note that we do not assume $d = \VC(\mathcal{H})$; in particular, it is possible that $d > \ell$. 

        \medskip

        \textbf{Base case:} $\ell \leq \ell_0$.

        \medskip

        We apply Theorem~\ref{theo:KKBell}, i.e., the Kahn--Kalai conjecture. If \ref{item:KK1} holds, then we are done as it implies \ref{item:VC1ind}. On the other hand, \ref{item:KK2} yields \ref{item:VC2ind} because
        \[
            48 q \log(\ell/\epsilon) \,\leq\, 48 q \log (300\, d^3/\epsilon^4) \,< \, C q \log(d/\epsilon),
        \]
        where in the last step we use that $C$ is large enough.

        \medskip

        \textbf{Inductive step:} $\ell > \ell_0$ and the theorem holds for all smaller values of $\ell$.

        \medskip

        Now $\ell > d$, so Lemma~\ref{lemma:count} can be applied. Pick $W_0 \sim X_{2q}$, set $t \defeq 3d \log \ell$, and let
        \[
            \mathcal{H}'\,\defeq\, \mathcal{H}_{W_0,t}^\mathsf{small} \qquad \text{and} \qquad \mathcal{F}_0 \,\defeq\, \mathcal{H}_{W_0,t}^\mathsf{large}.
        \]
        (This notation is defined in \S\ref{subsec:double}.) By Lemma~\ref{lemma:count} and Markov's inequality,
        \begin{equation}\label{markov}
           \P\left[\sum_{F \in \mathcal{F}_0} q^{|F|} \,\geq\, 2 \left(\frac{e}{d\ell}\right)^d\right] \,\leq\, \ell^{2d} \, 2^{-t} \,=\, \ell^{-d}. 
        \end{equation}
        By parts \ref{item:small} and \ref{item:VC} of Lemma~\ref{property}, $\mathcal{H}'$ is $t$-bounded and $\VC(\mathcal{H}') \leq d$. Since $\ell > \ell_0$, we have $t < \ell$, so, by the inductive hypothesis applied to $\mathcal{H}'$, at least one of the following holds:
        \begin{enumerate}[label=\ep{$\text{\normalfont{\small{Ind}}}'_{\arabic*}$}]
            \item\label{item:VC1prime} There is a set system $\mathcal{F}'$ with
            \[
                \mathcal{H}' \,\subseteq\, \upset{\mathcal{F}'} \quad \text{and} \quad \sum_{F \in \mathcal{F}'} q^{|F|} \,\leq\, \frac{1}{2} + 2\sum_{i = \ell_0}^t \left(\frac{e}{di}\right)^d.
            \]

            \item\label{item:VC2prime} For $p' = C q (\log(d/\epsilon) + \lambda_d(t))$ and $W' \sim X_{p'}$, we have 
            \[
                \P[W' \in \upset{\mathcal{H}'}] \,>\, 1 - \epsilon - \sum_{i = \ell_0}^t i^{-d}.
            \]
        \end{enumerate}
        We emphasize that whether \ref{item:VC1prime} or \ref{item:VC2prime} holds depends on the random set $W_0$.

        \medskip

        \textbf{Case 1:} \ref{item:VC1prime} happens with probability $\alpha > \ell^{-d}$.

        \medskip

        When \ref{item:VC1prime} holds, we let $\mathcal{F}'$ be a set system as in \ref{item:VC1prime} and define $\mathcal{F} \defeq \mathcal{F}_0 \cup \mathcal{F}'$. Then, by Lemma~\ref{property}\ref{item:undercover}, $\mathcal{H} \subseteq \upset{\mathcal{F}}$. Moreover, by \eqref{markov}, with probability at least $\alpha - \ell^{-d} > 0$,
        \[
            \sum_{F \in \mathcal{F}} q^{|F|} \,\leq\, \frac{1}{2} +2\sum_{i = \ell_0}^t \left(\frac{e}{di}\right)^d + 2 \left(\frac{e}{d\ell}\right)^{d} \,\leq\, \frac{1}{2} + 2\sum_{i = \ell_0}^\ell \left(\frac{e}{di}\right)^d,
        \]
        and hence \ref{item:VC1ind} is satisfied.

        \medskip

        \textbf{Case 2:} \ref{item:VC1prime} happens with probability at most $\ell^{-d}$.

        \medskip

        Set $p' \defeq C q (\log (d/\epsilon) + \lambda_d(t))$, choose $W' \sim X_{p'}$ independently of $W_0$, and let $W \defeq W_0 \cup W'$. Note that $W \sim X_{\tau}$ for some
        \begin{equation}\label{eq:tau_bound}
            \tau \,\leq\, 2q + p' \,=\, C q \log (d/\epsilon) + Cq\lambda_d(3 d \log \ell) + 2q. 
        \end{equation}
        We claim that, assuming $C$ is sufficiently large,
            \begin{equation}\label{eq:lambdas}
                \lambda_d(3 d \log \ell) + \frac{2}{C} \,\leq\, \lambda_d(\ell).
            \end{equation}
        Indeed, note that $\ell > \ell_0 > 9d^2$. If $\ell \leq 2^{3d}$, then $3d \log \ell \leq 9d^2$, which means that in this case
        \[
            \lambda_d(\ell) \,=\, \log^\ast(\ell) + 1 \,\geq\, \log^\ast(3d\log \ell) + 1 \,=\, \lambda_d(3d\log\ell) + 1,
        \]
        and we are done for $C \geq 2$. On the other hand, if $\ell > 2^{3d}$, then we use Claim~\ref{logrec} to write
        \begin{align*}
            \lambda_d(\ell) \,&=\, \log^\ast \ell + 2 \,\geq\, \log^\ast (\log^2 \ell) + 1/2 + 2 \\
            &\geq\, \log^\ast (3d \log \ell) + 1/2 + 2 \,\geq\, \lambda_d(3d \log \ell) + 1/2,
        \end{align*}
        and we are done assuming $C \geq 4$.
        
        Putting \eqref{eq:tau_bound} and \eqref{eq:lambdas} together, we see that
        \[
            \tau \,\leq\, Cq(\log (d/\epsilon) + \lambda_d(\ell)) \,=\, p.
        \]
        Now, using Lemma~\ref{property}\ref{item:induction}, we obtain
        \[
            \P[W \in \upset{\mathcal{H}}] \,\geq\, \P[\text{\ref{item:VC2prime}}] \,\P[W' \in \upset{\mathcal{H}'} \,\vert\, \text{\ref{item:VC2prime}}] \,\geq\, 1 - \ell^{-d} - \epsilon - \sum_{i = \ell_0}^t i^{-d} \,\geq\, 1 - \epsilon - \sum_{i=\ell_0}^\ell i^{-d},
        \]
        which yields \ref{item:VC2ind} and finishes the proof.
    \end{scproof}

    Now it is easy to derive Theorem~\ref{theo:KK}, which we restate here for convenience:

    \begin{theocopy}{theo:KK}
        There is a universal constant $A > 0$ with the following property. Let $\mathcal{H}$ be an $\ell$-bounded set system on a ground set $X$ with $\VC(\mathcal{H}) \leq d$ and let $q \in [0,1]$. Then at least one of the following statements holds.
        \begin{enumerate}[label=\ep{$\text{\normalfont{\small{VC}}}_{\arabic*}$}]
            \item\label{item:VC1again} There is a set system $\mathcal{F}$ with $\mathcal{H} \subseteq \upset{\mathcal{F}}$ and $\displaystyle \sum_{F \in \mathcal{F}} q^{|F|} \leq \frac{2}{3}$.

            \item\label{item:VC2again} For all $\epsilon \in (0, \frac{1}{2}]$, if $p = A q (\log(d/\epsilon) + \log^\ast \ell)$ and $W \sim X_p$, then $\P[W \in \upset{\mathcal{H}}] > 1 - \epsilon$.
        \end{enumerate}
    \end{theocopy}
    \begin{scproof}
        Without loss of generality, we may assume that $d \geq 2$. Let $C$ be the constant from Theorem~\ref{theo:KKvariant}. We will show that Theorem~\ref{theo:KK} holds with $A = 4C$. Apply Theorem~\ref{theo:KKvariant} with $\epsilon/2$ in place of $\epsilon$. Then $\ell_0=2400(d/\epsilon)^3$ and \ref{item:VC2ind} is stated for $W \sim X_{p'}$, where
        \[
            p' \,=\, Cq(\log(2d/\epsilon) + \lambda_d(\ell)) \,\leq\, Cq(\log(d/\epsilon) + \log^\ast(\ell) + 3) \,\leq\, p.
        \]
        Thus, to verify that \ref{item:VC1ind} and \ref{item:VC2ind} imply \ref{item:VC1again} and \ref{item:VC2again}, it is enough to check that
        \begin{align*}
            &2 \sum_{i = \ell_0}^\ell \left(\frac{e}{di}\right)^d \,\leq\, 2 \sum_{i = \ell_0}^\ell \left(\frac{e}{di}\right)^2 \,\le \, \frac{3e^2}{d^2\ell_0} 
            \,\le \, \frac{e^2\epsilon^3}{800d^5}
             \,<\, \frac{1}{6},\\
            &\sum_{i = \ell_0}^\ell i^{-d} \,\leq\, \sum_{i = \ell_0}^\ell i^{-2} \,< \, \frac{2}{\ell_0} \,<\, \frac{\epsilon^3}{1200 d^3} \,<\, \frac{\epsilon}{2}. \qedhere
        \end{align*}
    \end{scproof}

    \section{Sunflowers in low-dimensional families: Proof of Theorem~\ref{theo:sunflower}}\label{stheo}

    Let us state the theorem again:

    \begin{theocopy}{theo:sunflower}
        There exists a universal constant $C > 0$ such that the following holds. Let $\mathcal{H}$ be an $\ell$-bounded set system with $\VC(\mathcal{H}) \leq d$, where $d \geq 2$. If
        \[
            |\mathcal{H}|^{1/\ell} \,\geq\, C r (\log d + \log^\ast \ell),
        \]
        then $\mathcal{H}$ contains an $r$-sunflower.
    \end{theocopy}
    \begin{scproof}
        We fix $d$ and proceed by induction on $\ell$. The base case $\ell \leq 1$ is trivial. Suppose $\ell > 1$ and the statement holds for all smaller values. Let $C \defeq 10A$, where $A$ is the constant from Theorem~\ref{theo:KK}, and set
    \[
        Q \,\defeq\, C r (\log d + \log^\ast \ell), \text{ so } |\mathcal{H}| \,\geq\, Q^\ell.
    \]
    The constant $C$ is chosen so that $Q \geq 2A r (\log (2d) + \log^\ast \ell)$. Therefore, applying Theorem~\ref{theo:KK} with $p = 1/(2r)$ and $\epsilon = 1/2$, we arrive at two cases.

    \medskip

    \textbf{Case 1:} For $W \sim X_{1/(2r)}$ we have $\P[W \in \upset{\mathcal{H}}] > 1/2$.

    \medskip

    By partitioning $X$ into $2r$ subsets uniformly at random, we see that in expectation at least $r$ of them contain members of $\mathcal{H}$. This implies that $\mathcal{H}$ must contain $r$ pairwise disjoint sets (i.e., an $r$-sunflower with an empty kernel).

    \medskip

    \textbf{Case 2:} There is a set system $\mathcal{F}$ with
            \[
                \mathcal{H} \,\subseteq\, \upset{\mathcal{F}} \quad \text{and} \quad \sum_{F \in \mathcal{F}} q^{|F|} \,\leq\, \frac{2}{3},
            \] 
    where $q \defeq Q^{-1}$. Note that $\0 \not \in \mathcal{F}$ since $1 > 2/3$. For each $i \geq 1$, we let
    \[
        \mathcal{F}_i \,\defeq\, \set{F \in \mathcal{F} \,:\, |F| = i}.
    \]
    For each $F \in \mathcal{F}$, let $\mathcal{H}_F$ be the set of all $S \in \mathcal{H}$ such that $F \subseteq S$. If for some $F \in \mathcal{F}_i$, we have $|\mathcal{H}_F| \geq Q^{\ell - i}$, then we may apply the inductive hypothesis to the set system $\set{S \setminus F \,:\, S \in \mathcal{H}_F}$ to find an $r$-sunflower in $\mathcal{H}_F$ (and hence in $\mathcal{H}$) with a kernel including $F$. Thus, we may assume that $|\mathcal{H}_F| < Q^{\ell - i}$ for all $F \in \mathcal{F}_i$. Since $\mathcal{H} = \bigcup_{F \in \mathcal{F}} \mathcal{H}_F$, this yields
    \[
        Q^\ell \,\leq\, |\mathcal{H}| \,<\, \sum_{i = 1}^\infty Q^{\ell - i} |\mathcal{F}_i| \,=\, Q^\ell \sum_{i=1}^\infty q^i |\mathcal{F}_i| \,=\, Q^\ell \sum_{F \in \mathcal{F}} q^{|F|} \,\leq\, \frac{2}{3} \, Q^\ell,
    \]
    a contradiction.
    \end{scproof}

    \section{Kahn--Kalai for 1-dimensional families: Proof of Theorem~\ref{theo:KKdim1}}\label{sd1genr}

    Let us state Theorem~\ref{theo:KKdim1} again:

    \begin{theocopy}{theo:KKdim1}
        There is a universal constant $A > 0$ with the following property. Let $\mathcal{H}$ be a set system on a ground set $X$ with $\VC(\mathcal{H}) \leq 1$ and let $q \in [0,1]$. Then at least one of the following statements holds.
        \begin{enumerate}[label=\ep{$\text{\normalfont{\small{1d}}}_{\arabic*}$}]
            \item\label{item:VC1dim11} There is a set system $\mathcal{F}$ with $\mathcal{H} \subseteq \upset{\mathcal{F}}$ and $\displaystyle \sum_{F \in \mathcal{F}} q^{|F|} \leq \frac{2}{3}$.

            \item\label{item:VC2dim11} For all $\epsilon \in (0, \frac{1}{2}]$, if $p = A q \log(1/\epsilon)$ and $W \sim X_p$, then $\P[W \in \upset{\mathcal{H}}] > 1 - \epsilon$.
        \end{enumerate}
    \end{theocopy}

    Theorem~\ref{theo:KKdim1} is deduced from the following modification of Lemma~\ref{lemma:count} (see \S\ref{subsec:double} for the relevant definitions):

    \begin{lemma}\label{lemma:count1}
        Suppose $\VC(\mathcal{H}) \leq 1$ and for each pair of elements $x$, $y \in X$, there is some set $S \in \mathcal{H}$ with $S \cap \set{x,y} = \0$. Let $p$, $q \in [0,1]$ be such that $p \geq 2q$ and let $W \sim X_p$ be a random subset of $X$. Then, for all $t \geq 0$,
        \[
            \E\left[\sum_{F \,\in\, \mathcal{H}_{W,t}^{\mathsf{large}}} q^{|F|}\right] \,\leq\, 2 \, \left(\frac{q}{p}\right)^{t}.
        \]
    \end{lemma}

    Crucially, compared to Lemma~\ref{lemma:count}, there is no dependence on $\ell$ in Lemma~\ref{lemma:count1}. 

    The proof of Lemma~\ref{lemma:count1} is almost exactly the same as that of Lemma~\ref{lemma:count}. The only difference is in the step where we bound the cardinality of the set
    \[
        \left|\left\{(W, F) \,:\, |W| = w, \, F \in \mathcal{H}_{W,t}^{\mathsf{large}}, \, |F| = k\right\}\right|.
    \]
    As in the proof of Lemma~\ref{lemma:count}, we say that $|F \setminus W| = j$, where $t \leq j \leq k$, and then start by fixing the set $Z \defeq W \cup F$, for which there are at most ${n \choose w + j}$ choices. We now claim that the set $Z$ determines $F$ uniquely. Otherwise, there would be two distinct sets $F$, $F'$ of size $k$ such that $F = S \cap \core{Z}{\mathcal{H}}$ and $F' = S' \cap \core{Z}{\mathcal{H}}$ for some $S$, $S' \in \mathcal{H}$. Also, let $S'' \in \mathcal{H}$ be any set such that $S'' \subseteq Z$, which exists since $Z \in \upset{\mathcal{H}}$. By the definition of the $\mathcal{H}$-core, $\core{Z}{\mathcal{H}} \subseteq S''$. Since $|F| = |F'|$ and $F \neq F'$, we can pick elements $x \in F \setminus F'$ and $y \in F' \setminus F$. Now observe that $S \cap \set{x,y} = \set{x}$, $S' \cap \set{x,y} = \set{y}$, and $S'' \cap \set{x,y} = \set{x,y}$. Furthermore, by the assumption of the lemma, there is some $S^* \in \mathcal{H}$ with $S^* \cap \set{x,y} = \0$. This shows that $\VC(\mathcal{H}) \geq 2$, a contradiction. Thus, once we know $Z$, we also know $F$. Given $F$, we have at most ${k \choose j}$ choices for the set $F \setminus W$, which then determines $W$. To summarize,
    \[
        \left|\left\{(W, F) \,:\, |W| = w, \, F \in \mathcal{H}_{W,t}^{\mathsf{large}}, \, |F| = k\right\}\right| \,\leq\, \sum_{j = t}^k {n \choose w + j} {k \choose j}.
    \]
    The remainder of the proof of Lemma~\ref{lemma:count1} proceeds in exactly the same way as the proof of Lemma~\ref{lemma:count}, so we omit it.

    Now we show how to derive Theorem~\ref{theo:KKdim1} from Lemma~\ref{lemma:count1}.

    \begin{scproof}[ of Theorem~\ref{theo:KKdim1}]
        The argument is similar to the proof of Theorem~\ref{theo:KKvariant} (and even somewhat simpler as there is no induction). Let $A > 0$ be a constant that we will assume to be sufficiently large for all of the inequalities that appear in the sequel to hold. 
        
        Suppose that there are two elements $x$, $y \in X$ such that every set in $\mathcal{H}$ contains at least one of them, i.e., $\mathcal{H} \subseteq \upset{\set{\set{x},\set{y}}}$. If $q \leq 1/4$, then we can take $\mathcal{F} = \set{\set{x}, \set{y}}$ to make \ref{item:VC1dim11} hold. On the other hand, for $A > 4$, if $q > 1/4$, then $p > 1$, so in that case \ref{item:VC2dim11} holds. Therefore, we may assume that such elements $x$, $y$ do not exist.
        
        Now we  pick $W_0 \sim X_{2q}$, set $t \defeq 100/\epsilon$, and let
        \[
            \mathcal{H}'\,\defeq\, \mathcal{H}_{W_0,t}^\mathsf{small} \qquad \text{and} \qquad \mathcal{F}_0 \,\defeq\, \mathcal{H}_{W_0,t}^\mathsf{large}.
        \]
        By Lemma~\ref{lemma:count1} and Markov's inequality,
        \begin{equation}\label{markov1}
           \P\left[\sum_{F \in \mathcal{F}_0} q^{|F|} \,\geq\, 2^{1-t/2}\right] \,\leq\, 2^{-t/2}. 
        \end{equation}
        The family $\mathcal{H}'$ is $t$-bounded, so, by Theorem~\ref{theo:KKBell}, at least one of the following holds:
        \begin{enumerate}[label=\ep{$\text{\normalfont{\small{KK}}}'_{\arabic*}$}]
            \item\label{item:VC1prime1} There is a set system $\mathcal{F}'$ with $\mathcal{H}' \subseteq \upset{\mathcal{F}'}$ and $\displaystyle \sum_{F \in \mathcal{F}'} q^{|F|} \leq \frac{1}{2}$.

            \item\label{item:VC2prime1} For $p' = 48 q \log(2t/\epsilon)$ and $W' \sim X_{p'}$, we have $\displaystyle \P[W' \in \upset{\mathcal{H'}}] > 1 - \frac{\epsilon}{2}$.
        \end{enumerate}
        As in the proof of Theorem~\ref{theo:KKvariant}, which of \ref{item:VC1prime1} and \ref{item:VC2prime1} holds depends on the set $W_0$.

        \medskip

        \textbf{Case 1:} \ref{item:VC1prime1} happens with probability $\alpha > 2^{-t/2}$.

        \medskip

        When \ref{item:VC1prime1} holds, we let $\mathcal{F}'$ be a set system as in \ref{item:VC1prime1} and define $\mathcal{F} \defeq \mathcal{F}_0 \cup \mathcal{F}'$. Then $\mathcal{H} \subseteq \upset{\mathcal{F}}$ and, by \eqref{markov1}, with probability at least $\alpha - 2^{-t/2} > 0$,
        \[
            \sum_{F \in \mathcal{F}} q^{|F|} \,\leq\, \frac{1}{2} + 2^{1-t/2} \,<\, \frac{2}{3},
        \]
        and hence \ref{item:VC1dim11} is satisfied.

        \medskip

        \textbf{Case 2:} \ref{item:VC1prime1} happens with probability at most $2^{-t/2}$.

        \medskip

        Set $p' \defeq 48 q \log (2t/\epsilon)$, pick $W' \sim X_{p'}$ independently of $W_0$, and let $W \defeq W_0 \cup W'$. Note that $W \sim X_{\tau}$ for some
        \[
            \tau \,\leq\, 2q + p' \,=\, 48 q \log (200/\epsilon^2) + 2q \,\leq\, A q \log(1/\epsilon) \,=\, p, 
        \]
        assuming $A$ is large enough. Finally, using Lemma~\ref{property}\ref{item:induction}, we obtain
        \[
            \P[W \in \upset{\mathcal{H}}] \,\geq\, \P[\text{\ref{item:VC2prime1}}] \,\P[W' \in \upset{\mathcal{H}'} \,\vert\, \text{\ref{item:VC2prime1}}] \,\geq\, 1 - 2^{-t/2} - \frac{\epsilon}{2} \,\geq\, 1 - \epsilon,
        \]
        which yields \ref{item:VC2dim11} and finishes the proof.
    \end{scproof}
    
	\printbibliography
    
\end{document}